\documentclass{amsart}

\usepackage{amssymb}
\usepackage{latexsym}
\usepackage{amsmath}
\usepackage{enumerate}
\usepackage{amsthm}
\usepackage{pifont}
\usepackage{tikz}
\usepackage{upgreek}
\usepackage{hyperref}

\newcommand{\qee} {\hspace*{2mm}\hfill \ding{109}}

\renewcommand{\iff}{\leftrightarrow}
\renewcommand{\leq}{\leqslant}
\renewcommand{\geq}{\geqslant}
\renewcommand{\preceq}{\preccurlyeq}

\renewcommand{\phi}{\varphi}

\let\conOmega\Omega

\renewcommand{\Theta}{\varTheta}
\renewcommand{\Phi}{\varPhi}
\renewcommand{\Psi}{\varPsi}
\renewcommand{\Xi}{\varXi}
\renewcommand{\Omega}{\varOmega}
\renewcommand{\Gamma}{\varGamma}

\newtheorem{theorem}{Theorem}[section]
\newtheorem{define}[theorem]{Definition}

\newtheorem{exa}[theorem]{Example}

\newtheorem{exerc}[theorem]{Exercise}

\newtheorem{conj}[theorem]{Conjecture}

\newtheorem{ques}[theorem]{Open Question}

\newtheorem{lem}[theorem]{Lemma}

\newtheorem{cor}[theorem]{Corollary}

\newtheorem{rem}[theorem]{Remark}

\newcommand{\mc}[1]{\mathcal {#1}}
\newcommand{\num}[1]{{\underline {#1}}}
 \newcommand{\tupel}[1]{{\langle #1 \rangle}}
\newcommand{\verz}[1]{\{ #1 \}}
\newcommand{\gnum}[1]{{\ulcorner #1 \urcorner}}
\newcommand{\bles}{\mathbin{<}}
\newcommand{\bleq}{\mathbin{\leq}}
\newcommand{\mf}[1]{{\mathfrak {#1}}}
 \newcommand{\nrhd}{\mathrel{\not\! \rhd}}
 \newcommand{\bin}{\mathbin{\in}}
 \newcommand{\To}{\Rightarrow}

\newcommand{\kraj}{Kraj\'{\i}{\v{c}}ek}
\newcommand{\sot}{\ensuremath{{\sf S}^1_2}}
\newcommand{\ocol}{:}
\renewcommand{\descriptionlabel}[1]%
{\hspace{\labelsep}\emph{#1:}}

\title{Self-Interpretations and Kraj\'{\i}{\v{c}}ek-Theories}

\author{Albert Visser}
 \address{Philosophy, Faculty of Humanities,
                Utrecht University,
               Newton House, attn. Utrecht University,
3rd floor, Newtonlaan 201, 
3584 BH Utrecht, the Netherlands
}
\email{a.visser@uu.nl}
\date{\today}

\begin{document}



\thanks{We are grateful to the audience of the Interpretations workshop
Lugano, October 8-9, 2026 and to
Mateusz {\L}e{\l}yk for his talk at this workshop. The interest of the audience and Mateusz' presentation made me
think about interpretation categoricity and related concepts. I thank Mateusz {\L}e{\l}yk and Benjamin Zayton for their
interest in this project and their comments. Benjamin provided some helpful examples.}

\begin{abstract}
We study self-interpretations of theories and their connection to \kraj-theories.
Our two main results are the following. Firstly, if the self-inter\-preta\-tions of a theory form a
group in the category ${\sf INT}_3$ of theories and interpretations, then the theory is a
\kraj-theory. Secondly, a sequential theory is a \kraj-theory iff all its self-interpretations are 
\emph{initial}. Here  initiality of the interpretation $K$ means roughly that its internal
interpretations of the weak theory \sot\ are  initial in all interpretations of \sot.
\end{abstract}

\maketitle

\section{Introduction}
In this paper, we study \emph{self-interpretations} of theories and their connection with \emph{\kraj-theories}. 
Self-interpretations can best be seen as endomorphisms in categories of theories and interpretations.
In this paper, we focus on two such categories.
The first is ${\sf INT}_1$, the category in which
two interpretations $K$ and $M$ of $V$ in $U$ are the same iff  there is a $U$-definable, $U$-verifiable  isomorphism between
$K$ and $M$. The second is ${\sf INT}_3$, the category in which
two interpretations $K$ and $M$ of $V$ in $U$ are the same iff, for all $V$-sentences $\phi$, we have
$U \vdash \phi^K \iff \phi^M$. 
There are three properties of the self-interpretations of a theory that we will consider.
\begin{itemize}
\item
\emph{All self-interpretations are equal to the identity.} This means that the endomorphism monoid of
the given theory trivialises. A theory that has this property in ${\sf INT}_i$ will be called \emph{$\imath_i$-rigid}.
The $\imath_1$-rigid theories will also be called \emph{interpretation categorical}. 
The $\imath_3$-rigid theories will also be called \emph{sententially interpretation-rigid}.
\item
\emph{All self-interpretations are invertible.} In other words, the endomorphism monoid forms a group.
 A theory that has this property in ${\sf INT}_i$ will be called \emph{$\imath_i$-invertible}.
 We note that the invertible interpretations in ${\sf INT}_1$ are bi-interpretations and that the
 invertible interpretations in ${\sf INT}_3$ are sentential congruences (also: elementary congruences).
 \item
\emph{All self-interpretations are initial.} The property of initiality will be introduced in Section~\ref{mosterdsmurf}.
Roughly, $K:U \rhd U$ is  initial iff, for all interpretations of the weak theory \sot\ in $U$, there is an interpretation $M$ of
\sot\ in $U$ such that $K\circ M$ is below $N$.
We will show that a sequential theory is a \kraj-theory iff all its self-interpretations are initial.

We will also introduce a notion of \emph{strong initiality} and prove that a 
sequential theory is a \kraj-theory iff all its self-interpretations are strongly initial.
 \end{itemize}

\kraj-theories form a broad group.  The basic idea is, roughly,
 that a \kraj-theory proves an instance of its own restricted inconsistency for every interpretation of a designated
weak arithmetical theory. We note that, for example, all decidable theories are \kraj-theories, since they simply do
not interpret the designated weak arithmetical theory.
The \kraj-theories cannot be foundational theories like Peano Arithmetic
or Zermelo-Fraenkel Set Theory because of their $\Sigma^0_1$-unsoundness. 
Still, we submit, they are very interesting theories that contribute to getting a full picture
of what is possible in the presence of enough coding.

Our paper has two main results. First, in Section~\ref{rigtokraj}, 
we show that every $\imath_3$-invertible theory is a \kraj-theory. 
Then, in Section~\ref{mosterdsmurf}, we prove a characterisation of sequential \kraj-theories in terms of the initiality of their self-interpretations.
We also prove a slightly different version of the characterisation in terms of \emph{strong initiality}.

One reason to be interested in  invertibility is that it is implied by \emph{rigidity}.
We easily see that $\imath_1$-rigidity implies $\imath_3$-rigidity and that $\imath_3$-rigidity implies
$\imath_3$-invertibility.
The concept of  $\imath_1$-rigidity is a salient notion. It was recently studied by Mirko Engler and Benjamin Zayton in their paper
\cite{engl:rigi26}. Engler and Zayton provide many examples in addition to some discussion of the meaning of the concept.

An application of rigidity is the insight that, if a theory is $\imath_1$-rigid and interpretable in a finite sub-theory,
then it is finitely axiomatisable.\footnote{In fact, we already have:   if a theory is $\imath_1$-invertible and interpretable in a finite sub-theory,
then it is finitely axiomatisable.} See  \cite{zoet:inte15} and \cite{engl:rigi26}.
Fedor Pakhomov and Alexander Zapryagaev in \cite{pakh:mult20b} show that
Presburger Arithmetic is $\imath_1$-rigid. In combination with the fact that it is not finitely axiomatisable, we see that Presburger Arithmetic is not
interpretable in a finite sub-theory. See \cite[III.4]{smor:logi91} (especially Exercise 13)
 and \cite{zoet:inte15} for a proof of non-finite axiomatisability. In contrast, Presburger 
Arithmetic is locally interpretable in a finitely axiomatised sub-theory.\footnote{Peano Arithmetic is not even locally interpretable in a
finitely axiomatised sub-theory.}

\subsection{Contents}
We introduce \kraj-theories in Section~\ref{vogelsmurf}. In this section we also present some salient earlier results.
For example, we repeat a result from \cite[Appendix C]{viss:faith05} that shows that the predicate logic of a c.e. \kraj-theory is
$\Pi^0_2$-complete. 
In Section~\ref{rigtokraj}, we prove that $\imath_3$-invertibility implies being a \kraj-theory. Finally, in Section~\ref{mosterdsmurf},
we characterise sequential \kraj-theories in terms of the behaviour of their self-interpretations.

\subsection{AI Statement}
ChatGPT and Claude diligently found typos and small mistakes in previous versions of this paper.
Also, they gave helpful naming advice. Chat GPT  provided
Example~\ref{rechtlijnigesmurf} and it directed my attention to Theorem~\ref{abstractesmurf}.
All the further results of the paper are my own.

\section{Basic Facts and Basic Definitions}

In this section, we lay the groundwork for the rest of the paper. Of course, the reader can touch upon each subsection
lightly, to return to it when needed. In Section~\ref{mosterdsmurf}, we will need  more heavy machinery.
We will introduce these further instruments in the section itself.

\subsection{Theories}
A theory in this paper will be a c.e. theory in predicate logic of finite signature. 
We will assume that the axiom sets of theories are simple, say given by a $\Delta^{\sf b}_1$-predicate, so
that we do not have to worry that the representation of the axiom set adds extra complexity to the
proof-predicate.

The theories {\sf R} and \sot\ will play a special role in this paper. The theory \sot\ was introduced by Sam Buss in \cite{buss:boun86}. See also \cite{haje:meta91}.
Buss' theory is a weak theory that is mutually interpretable with such theories as {\sf Q} and ${\sf PA}^-$. On the other hand, it is an ideal basis for the arithmetisation
of metamathematics. The theory {\sf R} was introduced in \cite{tars:unde53}. Since, we will see some actual reasoning in {\sf R}, we reproduce its axioms here.
For $m, n \in \upomega$ we have: 
\begin{enumerate}[{\sf R}1.]
	\item  $\num{m} + \num{n} = \num{m + n}$
	\item  $\num{m} \times \num{n} = \num{m \times n}$
	\item  $\num{m} \neq \num{n}$ (if $m \neq n$)
	\item  $\forall x\, \bigl(x \leq \num{n} \to \bigvee_{i \leq n} x = \num{i}\bigr)$
	\item  $\forall x\, (x \leq \num{n} \lor \num{n} \leq x)$
\end{enumerate}
In the Tarski-Mostowski-Robinson book, the relation $\leq$ is defined, but this is not really needed and we take it as primitive here.

\subsection{Interpretations and Interpretability}
Most of my papers contain a section describing the basics of interpretations and interpretability.
See, e.g., \cite[Appendix A]{viss:onq17}. We just give some highlights here.

An interpretation $K$ is based on a translation $\tau$ of languages. 

We write $K:U \rhd V$ for: $K$ is an interpretation of $V$ in $U$, and $U \rhd V$ for: there is an interpretation
$K:U \rhd V$. 

We will diverge a bit from the usual practice by also using $\tau: U \rhd V$ with as meaning: there is an
interpretation $K:U \rhd V$ with underlying translation $K$. The reason we use this is to avoid having to
name different interpretations with the same underlying translation.

We remind the reader of 
disjunctive interpretations. 
Suppose $K\ocol (U+\phi) \rhd V$ and $M\ocol (U+\neg\, \phi) \rhd V$. Then $K\tupel{\phi}M\ocol  U \rhd V$ is the
disjunctive interpretation that coincides with $K$ if $\phi$ and with $M$ if $\neg \,\phi$.

Suppose $U$ and $V$ are theories in the same language and the theorems of $V$ are contained in the theorems of $U$.
We write ${\sf E}_{VU}:U \rhd V$ for the interpretation based on the identical translation.

Suppose $N\ocol U \rhd {\sf R}$. An interpretation $K\ocol U \rhd V$ is \emph{$N$-restricted} iff there is a $U$-formula $\theta(x)$, such that,
for all $\phi$, we have $U \vdash \theta(\gnum \phi) \iff \phi^K$. Here the numerals come from $N$. 
The interpretation $K$ is \emph{restricted} iff it is $N$-restricted, for some $N$.
The present notion is  weaker than the notion of bounded interpretation from \cite{viss:seco12}.
We provide a stronger version in Section~\ref{mosterdsmurf}.
We write $K\ocol U \rhd_{\sf res} V$, for restricted interpretability.
We have:
\begin{theorem}
\begin{itemize}
\item
If $U \rhd V \rhd_{\sf res}W$, then $U \rhd_{\sf res}W$.
\item
If $U \rhd_{\sf res} V \rhd W$, then $U \rhd_{\sf res}W$.
\end{itemize}
\end{theorem}

We leave the simple verification to the reader.
Let
\begin{itemize}
\item
$\mho(U) := \sot+ \verz{{\sf con}_{\num n}(U)\mid n\in \upomega}$. 
\end{itemize}
Here ${\sf con}_n$ means consistency from the axioms with G\"odel numbers $\leq n$ where only formulas of
depth-of-quantifier complexity $\leq n$ are allowed in the proofs. We have:
\begin{theorem}[Interpretation Existence Lemma]\label{intexists}
We have: $\mho(U) \rhd_{\sf res} U$. 
\end{theorem}
For a proof see \cite{viss:inte18}. See also \cite{viss:seco12}. The present formulation is a convenient version for most of our
purposes,
however, we can articulate much stronger versions. We will meet a stronger version in Section~\ref{mosterdsmurf}, to wit Theorem~\ref{sterkesmurf}.

An interpretation is \emph{direct} if it is one-dimensional, its domain is the full domain of the interpreting theory and
it translates identity to identity of the interpreting theory. 

There are \emph{two categories of theories and interpretations} relevant for this paper.
In the context of the categories, we usually write $K\ocol V\to U$ for $K\ocol U \rhd V$.
The first category is ${\sf INT}_1$. In this category, two interpretations $K,M\ocol V \to U$ are the same if
there is an $U$-definable, $U$-verifiable isomorphism between $K$ and $M$. Isomorphism of theories in this
category is \emph{bi-interpretation}. 
The second is ${\sf INT}_3$. In this category, two interpretations $K,M\ocol V \to U$ are the same if, for all $V$-sentences $\psi$,
we have $U \vdash \psi^K \iff \psi^M$. 
 Isomorphism of theories in this
category is \emph{elementary congruence} or \emph{sentential congruence}.

\subsection{Sequentiality}
A theory is \emph{sequential} if it has a good theory of sequences, where sequences are sequences of all objects and
where we can internally add an element to a sequence. Via a long bootstrap, we can prove that
a theory is sequential iff it directly interprets the very weak set theory {\sf AS}. 
This was first proved in its most general form in \cite{myci:latt90}. 
See \cite{viss:what13} for an exposition and more of the history of the result.

An important property of sequential theories is \emph{local reflexivity}. Suppose $U$ is sequential. Then, for any $n$, there is
an $N\ocol U \rhd (\sot+ {\sf con}_n(U))$.
This property can be proven using the methodology of shortening cuts in combination with the existence of partial satisfaction predicates.
See, e.g., \cite[Section 2]{viss:smal19}.

\subsection{Witness Comparison}
Suppose $\phi$ is of the form $\exists x\, \phi_0(x)$ and $\psi$ is of the form $\exists y\, \psi_0(y)$.
We define:
\begin{itemize}
\item
$\phi \leq \psi :\iff \exists x\, (\phi_0(x) \wedge \forall y\bles x\, \neg\, \psi_0(y))$.
\item
$\phi < \psi :\iff \exists x\, (\phi_0(x) \wedge \forall y\bleq x\, \neg\, \psi_0(y))$.
\item
$(\phi \leq \psi)^\bot :\iff \psi < \phi$.
\item
$(\phi < \psi)^\bot :\iff \psi \leq \phi$.
\end{itemize}

\section{\kraj-Theories}\label{vogelsmurf}
The c.e. theory
$U$ is a \emph{\kraj-theory} iff, for every $N\ocol U \rhd {\sf S}^1_2$, there is a $k$ such that
$U \vdash {\sf incon}^N_k(U)$.
These theories are named after Jan \kraj, who introduced the basic idea in \cite[Section 3]{kraj:note87}.
Our notion is broader than the original one. First, we do not restrict ourselves to sequential theories and, secondly,
we do not restrict ourselves to theories constructed from finitely axiomatised ones (as is done in our Section~\ref{uitvoerigesmurf}).

One of \kraj's results in \cite{kraj:note87} implies that consistent sequential \kraj-theories cannot be
finitely axiomatisable.

We prove two basic theorems about the preservation of \kraj-theories.
\begin{theorem}
Suppose $U$ is a \kraj-theory and  $V$ is an ${\sf INT}_3$-retract of $U$. Then, $V$ is a \kraj-theory.
\end{theorem}

\begin{proof}
Let $K\ocol  U \rhd V$ and $M\ocol V \rhd U$ and, for all $V$-sentences $\psi$, $V \vdash \psi \iff \psi^{MK}$.
Suppose $N\ocol  V \rhd \sot$. It follows that $K \circ N\ocol U \rhd \sot$. So, for some $k$, we have
$U \vdash {\sf incon}_k^{NK}(U)$. Thus, $V \vdash {\sf incon}^{NKM}_{k}(U)$. 
We may conclude  $V \vdash {\sf incon}^N_{k}(U)$. 

Using $M$ we can transform a $U$-proof of inconsistency into a $V$-proof of inconsistency.
We note that the translated proof from the translated $U$-axioms can be executed in p-time.
This leads to a verification in \sot. Moreover, the translation only adds a standard constant to the complexity
of the formulas. Finally, we have only to worry about standardly many $V$-axioms. It follows that
  $V \vdash {\sf incon}^{N}_{k^\ast}(V)$, for some standard $k^\ast$. 
\end{proof}

The following result is a bit messy, but it does deliver the desired goods.

\begin{theorem}\label{minismurf}
Suppose $U$ is a sub-theory of $V$ \textup(qua theorems, in the same language\textup) and $U \rhd \sot$.
Suppose further that for every $N:U \rhd \sot$, there is an $n$ such that $V \vdash {\sf incon}^N_n(U)$.
Then, $V$ is a \kraj-theory.
\end{theorem}

We will use translations rather than interpretations to realise the proof.

\begin{proof}
We assume the conditions of the theorem.

Let $\nu^\ast: U \rhd \sot$. Suppose $\nu:V \rhd \sot$. Then, for some $\phi$ with $V \vdash \phi$, we have
$\nu:(U+\phi)\rhd \sot$. Let $\nu^\circ:= \nu\tupel{\phi}\nu^\ast: U \rhd \sot$. It follows that, for some $n$,
$V \vdash {\sf incon}^{\nu^\circ}_n(U)$ and, thus, $V \vdash {\sf incon}^{\nu}_n(U)$.
Suppose the $V$-axioms needed to prove the $U$-axioms $\leq n$ are all $\leq m$ and let $k := {\sf max}(m,n)$.
Then, $V \vdash {\sf incon}^{\nu}_k(V)$.
\end{proof}

\begin{cor}
Let $U$ be a sub-theory of $V$ and suppose $U$ is \kraj. Then, $V$ is \kraj.
\end{cor}

\subsection{Extension to a \kraj-Theory}\label{uitvoerigesmurf}

In this subsection we show that if $A$ is finitely axiomatized, then we can locally interpret a 
 \emph{Kraj\'{\i}{\v{c}}ek theory} extending  $A$ in $A$. This result, applied to sequential theories, establishes the existence of
 \kraj-theories in the sequential realm.
 
 We  first need an auxiliary result
 from \cite{viss:fefe14}, to wit a restricted version of Feferman's Theorem on the interpretability of
inconsistency for finitely axiomatised  theories. Let $A$ be  finitely axiomised.
We write $\uprho(\phi)$ for the depth-of-quantifier-alternations complexity of $\phi$ and $\uprho(N)$ for the maximum
of the $\rho$-complexities of the formulas involved in the indwelling translation of $N$. The standard constant $\mf c$ 
just depends on details of the coding and the direct interpretation witnessing sequentiality. This tells us that the result will be
preserved to finite extensions of $A$.

\begin{theorem}\label{salsaparilla}
Suppose $N:A \rhd {\sf S}^1_2$ and $n\geq {\sf max}(\uprho(A),\uprho(N)+\mathfrak c)$.
Then, we have $A \rhd (A+ {\sf incon}^N_{n}(A))$. 
\end{theorem}

With Theorem~\ref{salsaparilla} in hand, we prove the following theorem.

 \begin{theorem}\label{huismansmurf}
Let $A$ be finitely axiomatized. 
Then, a c.e. theory $U:= A+\verz{{\sf incon}^{N_k}_{n_k}(A) \mid k \in \upomega}$ is locally interpretable in $A$.
Here $N_k$ enumerates the $N:A \rhd \sot$ and $n_k$ is some effective sequence of natural numbers.
Moreover, $U$ is a \kraj-theory.
\end{theorem}

We note that we prove a bit more than \kraj\ in his \cite[Theorem 3.0]{kraj:note87}, since we do not
have to restrict ourselves to the sequential case. The reason for that is that Feferman's idea of
interpretability of inconsistency does not demand sequentiality.

\begin{proof} 
In case $A \nrhd {\sf S}^1_2$, the desired result is trivial. We take $U := A$.
Suppose $A \rhd {\sf S}^1_2$.
 Let  $N_k: A\rhd \sot$ enumerate the interpretations $N:A\rhd \sot$. The indwelling translations of the $N_k$ are $\nu_k$.
 Then, we can find sequences $K_0,K_1,\ldots$ and  $n_0,n_1, \ldots$ such that, for every
 $k$, we have:  $K_k:A \rhd (A +\bigwedge_{i<k} {\sf incon}_{n_i}^{N_i} (A))$. 
 The indwelling translations of the $K_k$ will be $\kappa_k$.
 
 We produce the $K_k$ and $n_k$ in stages. At stage $k$ we produce 
$K_k$ and at stage $k+1$ we construct $n_k$.
 
At stage 0, we take $K_0 := {\sf ID}_A$. We consider stage $k+1$.
We are given $K_k: A \rhd (A +\bigwedge_{i<k} {\sf incon}_{n_i}^{N_i} (A))$. 
Clearly, $ \kappa_k\circ \nu_k:A \rhd {\sf S}^1_2 $. By Theorem~\ref{salsaparilla}, we can effectively
find an $n_k$ and an $M$ such that $M: A \rhd (A + {\sf incon}_{n_k}^{\nu_k\kappa_k} (A))$.
Let the indwelling translation of $M$ be $\mu$.
We find that $A  \vdash  (\bigwedge_{i<k} {\sf incon}_{n_i}^{\nu_i}(A))^{\kappa_k}$.
So, it follows that:
\begin{eqnarray*}
 A & \vdash &  
 \big (\bigwedge_{i<k} {\sf incon}_{n_i}^{\nu_i}(A)\big )^{\kappa_k\mu} \wedge {\sf incon}_{n_k}^{\nu_k\kappa_k\mu} (A) \\
 & \vdash & \big(\bigwedge_{i<k+1} {\sf incon}_{n_i}^{\nu_i}(A)\big)^{\kappa_k\mu}
\end{eqnarray*}

\noindent
So, we find the desired $K_{k+1}$ with translation $\mu \circ \kappa_k$.
We may conclude that $A \rhd_{\sf loc} (A + \verz{{\sf incon}_{n_k}^{N_k}(A) \mid k\in \upomega})$.

Finally,  $U$ is indeed a \kraj-theory. If $A \nrhd \sot$, we are immediately done. 
If $A\rhd \sot$, the desired result follows from Theorem~\ref{minismurf}.
\end{proof}

 \noindent
 We note that if $\tau_i$ is an enumeration of translations of the signature or arithmetic to the signature of $A$
 and if $\sigma$ carries an interpretation of ${\sf S}^1_2$ in $A$, then we can effectively enumerate 
 translations supporting the $N_i$ as $\tau_i\tupel{({\sf S}^1_2)^{\tau_i}} \sigma$. Thus our axiomatization
 of $U$ can be made decidable. 
 
 \begin{rem}
 We note that every finitely axiomatised $A$ model-interprets in a sense a \kraj-variant.
 Consider any model $\mc M$ of $A$. In case all interpretations of $\sot$ in $\mc M$  contain 
 ${\sf incon}_{n}(A)$, for some $n$, the identical interpretation gives us a \kraj-extension of
 $A$. If not, there is an $N$ in $\mc M$ in which we have $\mho_A$. Using this we can construct an internal
 model of a \kraj-extension since local interpretability in $A$ implies global interpretability in $\mho_A$.
 However, in the first case, we do not always have the same theory. Moreover, these theories could fail to be c.e. 
 \end{rem}

\subsection{The Predicate Logic of a \kraj-Theory}

We consider the predicate logical schemes valid for all substitutions in $U$. We call these schemes $\Lambda_U$, where we suppress the
chosen signature. We could look at schemes for relativised interpretations of predicate logic but also at schemes
for direct interpretations, where we do not use the original theory $\sot$ but its translation to a relativised format where
identity becomes an equivalence relation. We will assume relativised interpretation here.

Let $\alpha$ be a single sentence expressing the axioms of \sot\ including the axioms of identity.
We can also use e.g. Adjunctive Set Theory, a theory with just one binary predicate.

\begin{theorem}\label{kleinesmurf}
Let $U$ be a consistent \kraj-theory. Suppose $U \rhd \sot$. Then, $\Lambda_U \vdash \alpha \to {\sf incon}(U)$.
\end{theorem}

\begin{proof}
Suppose $N_0\ocol U \rhd \sot$. Let $K$ be any interpretation of predicate logic in the desired signature in $U$.
Let $\tau_0$ be the indwelling translation of $N_0$ and let $\tau$ be the indwelling translation of $K$.
Consider $\tau^\ast := \tau\tupel{\alpha^\tau} \tau_0$. Clearly, $\tau^\ast$ supports an interpretation $K^\ast$ of \sot.
It follows that $U \vdash {\sf incon}^{K^\ast}(U)$, since $U \vdash {\sf incon}_k^{K^\ast}(U)$, for some $k$.
So, $U \vdash (\alpha \to  {\sf incon}(U))^K$.
\end{proof}

\begin{exa}
Consider predicate logic with at least a binary predicate symbol. Clearly, this is a \kraj-theory. However, its
schematic logic is simply ordinary predicate logic.
\end{exa}

The following result is taken from \cite[Appendix C]{viss:faith05}. We repeat it, since the results
of the present paper give it new relevance.

\begin{theorem}
Suppose $U$ is a consistent sequential \kraj-theory.
Then, $\Lambda_U$ is $\Pi^0_2$-complete.
\end{theorem}

\begin{proof}
Since $U$ is sequential, there is an $N_0\ocol U \rhd \sot$. We fix $N_0$ with indwelling translation $\tau_0$.

Let $\sigma(x)$ be any $\Sigma^0_1$-formula. We find a $\rho(x)$ such that 
\[\sot \vdash \rho(x) \iff \big (\sigma(x) \leq \exists N \,{\sf prov}_U(\gnum{(\alpha \wedge \rho(\dot x))^N})\big ).\]
Here we assume that the $ \exists N$ and the witness of {\sf prov} are combined in one existential quantifier.
We show that:
\begin{enumerate}[A.]
\item
$\sigma(\num n)$ is true iff $\rho(\num n)$ is true.
\item
 $\sigma(\underline n)$ is true iff, for some $N\ocol U \rhd \sot$,  we have  $U\vdash \rho^N(\underline n)$.
 \end{enumerate}
We first prove (A). The right-to-left direction is trivial. We show left-to-right.
Suppose $\sigma(\underline n)$ is true. Then,  $\rho(\underline n)$ is true or $\rho^\bot(\underline n)$ is true.
In the first case we are done. In the second case we have, $U\vdash \rho^N(\num n)$, for some $N\ocol U\rhd \sot$, by the definition of $\rho^\bot$.
But, by $\Sigma^0_1$-completeness, we also have $U\vdash \rho^{\bot N}(\num n)$. So, $U$ would be inconsistent.
\emph{Quod non.}

We prove (B). Suppose $\sigma(\num n)$ is true. Then, by (A), we find that $\rho(\num n)$ is true. So, $U \vdash \rho^{N_0}(\num n)$.
In the other direction, suppose  $N\ocol U \rhd \sot$,  and  $U\vdash \rho^N(\underline n)$. It follows that
$\rho(\num n)$ is true or $\rho^\bot(\num n)$ is true. In the first case we are immediately done. In the second case,
we get a contradiction via $\Sigma^0_1$-completeness as before.\footnote{The above argument is a minor variant of the Friedman-Goldfarb-Harrington
argument. The true priority for the argument goes to John Shepherdson in \cite{shep:rep61}.}

\medskip
We claim that 
\[\forall x\, \sigma(x)\text{ is true iff }\Lambda_U \vdash \alpha \to \exists y\, \big({\sf incon}_{y}(U) \wedge \forall x\bleq y\, \rho(x) \big).\]

Suppose that $\forall x\, \sigma(x)$ is true.  Let $K$ be any interpretation of predicate logic in the desired signature in $U$ with indwelling translation $\tau$.
We consider the translation $\tau^\ast: \tau \tupel{\alpha^\tau} \tau_0$. Clearly, $\tau^\ast$ supports an interpretation $K^\ast\ocol U \rhd \sot$.
It follows that, for some $k$, we have $U\vdash {\sf incon}^{K^\ast}_k(U) $. 
By (A), we have $\forall x\, \rho(x)$ is true. So, certainly: $U\vdash ({\sf incon}_k(U) \wedge \forall x \bleq \num k\, \rho(x))^{K^\ast}$ and, thus,
$U\vdash (\exists y\, ({\sf incon}_y(U) \wedge \forall x \bleq y\, \rho(x)))^{K^\ast}$. 
We may conclude that
$U \vdash \big(\alpha \to \exists y\, ({\sf incon}_y(U) \wedge \forall x \bleq y\, \rho(x))\big)^{K}$.

We prove the converse. Suppose (\dag) $\Lambda_U \vdash \alpha \to \exists y\, \big({\sf incon}_{y}(U) \wedge \forall x\bleq y\, \rho(x) \big)$.
We show that, for any $k$, the sentence $\forall x \bleq \num k\, \sigma(x)$ is true.
By local reflexivity, there is an $N\ocol U \rhd \sot$, such that $U \vdash {\sf con}_k^N(U)$. By (\dag), we find that
 \[U \vdash \big(\exists y\, ({\sf incon}_y(U) \wedge\forall x\bleq y\, \rho(x)) \big)^N .\] 
 Since $U \vdash (({\sf con}_{\num k}(U) \wedge {\sf incon}_y(U)) \to y<\num k)^N$,
 it follows that,
 $U \vdash (\forall x\bleq \num k\, \rho(x))^N$. By (B), we may conclude that $\forall x \bleq \num k\; \sigma(x)$ is true.
\end{proof}

We note that the above proof already works for locally reflexive \kraj-theories.

\section{Invertible Theories are \kraj}\label{rigtokraj}
In this section, we prove one of the two main results of this paper: invertible theories are \kraj.

\subsection{Basic Definitions}
We repeat the basic definitions.
\begin{itemize}
\item
The theory $U$ is $\imath_i$-rigid iff its endomorphism monoid in ${\sf INT}_i$ is trivial, i.e., all self-interpretations of
$U$ are ${\sf INT}_i$-equal to the identity. We say that $U$ is
 \emph{interpretation-categorical} if it is $\imath_1$-rigid and that it is \emph{sententially interpretation-rigid} if it is $\imath_3$-rigid.
 \item
 The theory $U$ is $\imath_i$-invertible iff its endomorphism monoid in ${\sf INT}_i$ is a group.
  We note that $U$ is $\imath_1$-invertible iff all its self-interpretations are bi-interpretations and that it
   is $\imath_3$-invertible iff all its self-interpretations are sentential congruences.
 \end{itemize}

The following example is a modification suggested by ChatGPT of a counter-example to rigidity proposed by Benjamin Zayton.

\begin{exa}\label{rechtlijnigesmurf}
Consider the theory {\sf DLO} of dense linear orderings without end-points. This theory is complete and hence decidable. So it is, trivially, 
$\imath_3$-rigid and also a \kraj-theory.

Let $K:{\sf DLO} \rhd {\sf DLO}$ be given by:
\begin{itemize}
\item
  $\updelta_K(x,y) :\iff \top$,
  \item
   $(x,y)=_K(u,v) :\iff (x=u \wedge y=v)$,
   \item
    $(x,y) \leq_K (u,v): \iff (x<u \vee (x=u \wedge y \leq v))$.
    \end{itemize}
Consider the models $\mathbb Q$ and $\mathbb Q\times_{\sf lex} \mathbb Q$. The interpretation $K$ yields a mapping $F:f \mapsto
\uplambda xy. (f(x),f(y))$ of automorphisms. It is easy to see that $F$ is not surjective. However, if $K$ had an ${\sf INT}_1$-inverse,
$F$ would have to be an isomorphism. So, {\sf DLO} is not $\imath_1$-invertible.
\end{exa}

The following basic result was pointed out to me by ChatGPT. 

\begin{theorem}\label{abstractesmurf}
Let $\mc C$ be a category and suppose the endomorphism monoid of $b$ is a group.
Suppose $a$ is a retract of $b$. Then, $a$ is isomorphic to $b$ and, hence, the endomorphism monoid of $a$ is also a group.
\end{theorem}

\begin{proof}
Suppose $f:a\to b$ and $g:b\to a$ and $g\circ f = {\sf id}_a$. It follows that $(f\circ g) \circ (f\circ g) =
f\circ (g \circ f) \circ g = f\circ g$. But, by assumption, $f\circ g$ is invertible, so it must be ${\sf id}_b$.
\end{proof}

Thus we may conclude that properties like $\imath_i$-rigidity and $\imath_i$-invertibility are preserved to
retracts, since the retracts will be isomorphisms.

\subsection{The Main Development}
We need two theorems as preparation for the main result of this section.

\begin{theorem}\label{slimmesmurf}
Let $U$ be a c.e. theory and let $N\ocol U \rhd {\sf R}$.
Suppose that the theory $U^\star :=U+ \verz{\phi(\underline n) \mid n \in \upomega}$ is consistent.
Here the $\underline n$ are the numerals provided by $N$.\footnote{The interpretation $N$ might be more dimensional
or piece-wise. Notationally, we pretend it is one-dimensional.} Then, there is sentence $\zeta$ such that
$U+\zeta$ is consistent and $U+\zeta \vdash  \verz{\phi(\underline n) \mid n \in \upomega}$.
\end{theorem}

\begin{proof}
Let $I(x) :\iff  \updelta_N(x) \wedge \phi(x)$. 
By the Fixed Point Lemma, we find a $\rho$ with 
\[{\sf R} \vdash \rho \iff \big ({\sf prov}_{U^\star}(\gnum{\neg\, \rho^I}) \leq {\sf prov}_{U^\star}(\gnum{\rho^I})\big ) .\] 
We can assume that both the substitution function and  the proof-predicate are elementary.
Moreover, the graph of exponentiation can be represented by a $\Delta_0$-formula in such a way that {\sf R} proves uniqueness
of values on standard arguments. It follows that $\rho$ can be written as $\exists x\, \rho_0(x)$, where $\rho_0$
is elementary. We find: if $\rho_0(\underline n)$ is true, then ${\sf R} \vdash  \rho_0(\underline n)$, and
if $\rho_0(\underline n)$ is false, then ${\sf R} \vdash \neg\, \rho_0(\underline n)$.

Our first order of business is to prove that $\rho^I$ is consistent with $U$ and that $\rho$  is false.

Suppose $U^\star \vdash \neg\, \rho^I$. Then, (i) $\rho$ is true or (ii) $\rho^\bot$ is true. 
\begin{enumerate}[\emph{C{a}se} i.]
\item
We find, by the usual proof
of $\Sigma^0_1$-completeness and the fact that the internal witness is standard, that $U^\star \vdash \rho^I$. This contradicts the fact
that $U^\star$ is consistent. 
\item
We have  $U^\star \vdash \rho^I$, by the definition of $\rho^\bot$. This again contradicts
the fact that $U^\star$ is consistent. 
\end{enumerate}
So, $U^\star$ is consistent with $\rho^I$ and, \emph{a fortiori}, $U$ is consistent with $\rho^I$.

If $\rho$ were true, we would have $U^\star \vdash \neg\, \rho^I$. We have already seen that this does not hold.

Let $\zeta := \rho^I$.  
We show that  $U+  \zeta \vdash  \verz{\phi(\underline n) \mid n \in \upomega}$.
Since $\rho$ is false, we find (\dag) $U \vdash \neg\, \rho^N_0(\underline n)$, for all $n$.
We work in $U+\zeta$. Let $z$ be a witness of $\rho^I$. Suppose $\neg\,\phi(\underline n)$. Then $\underline n \not\in I$.
So, $\underline n \not\leq z$ and, thus, by Axiom {\sf R}5, $z \leq \underline n$. It follows, by Axiom {\sf R}4,
  that $\bigvee_{k\leq n} z= \underline k$, and, thus, 
 $\bigvee_{k\leq n} \rho_0^N(\underline k)$, but this contradicts (\dag). So, $\phi(\underline n)$.
\end{proof}

We note that $\rho^I$ can also be written as $\rho^N< \exists y\bin \updelta_N \,\neg\, \phi(y)$.

\begin{rem}
We are not sure how original Theorem~\ref{slimmesmurf} is. Let's say that results in the same spirit are known.
\end{rem}

In the following theorem, the correct place where the numerals live is crucial, so we will keep track of it by adding a subscript.

\begin{theorem}\label{sluwesmurf}
Consider a c.e. theory $U$ and let $N\ocol U \rhd {\sf R}$. Suppose $U^\star := U+ \verz{\phi(\underline n^N) \mid n \in \upomega}$ is consistent and
$U^\star \rhd_{{\sf res}} U$.  Then, $U$ is not $\imath_3$-invertible.
\end{theorem}

\begin{proof}
By Theorem~\ref{slimmesmurf}, there is a $\zeta$ such that $U+\zeta$ is consistent
and $U+\zeta$ extends $U^\star$. Clearly, we have $(U+\zeta) \rhd_{\sf res} U$.
Let $K\ocol (U+\zeta) \rhd U$, $M\ocol (U+\zeta)\rhd {\sf R}$ and $\theta(x)$ witness the restricted interpretation.

Let $P:= K \tupel{\zeta}{\sf E}_{U(U+\neg\, \zeta)}\ocol U \rhd U$. 
Suppose $U$ is $\imath_3$-invertible. Let $Q:U \rhd U$ be the promised ${\sf INT}_3$-inverse.
We take  $M' :=Q\circ (M\tupel{\zeta}N):U \rhd {\sf R}$.
By the Fixed Point Lemma, we find $\lambda$ with $U \vdash \lambda \iff \neg \, \theta^Q(\gnum \lambda^{M'})$.
We have:
\begin{eqnarray*} 
U+\zeta^Q  \vdash \theta^Q(\gnum{\lambda}^{M'}) & \iff & \theta^Q(\gnum{\lambda}^{Q\circ M})  \\
& \iff &   \big(\theta(\gnum{\lambda}^{M})\big)^Q\\
& \iff & \lambda^{KQ} \\
& \iff & \lambda^{PQ} \\
& \iff & \lambda
\end{eqnarray*} 
So, $U + \zeta^Q \vdash \lambda \iff \neg\, \lambda$ and, hence, $U \vdash (\neg\,\zeta)^Q$.
On the other hand, $Q$ as a sentential congruence is faithful, so $U \vdash \neg\, \zeta$.
But this contradicts the consistency of $U+\zeta$.
\end{proof}

We prove our main result.

\begin{theorem}
Suppose $U$ is $\imath_3$-invertible. Then, $U$ is a \kraj-theory.
\end{theorem}

\begin{proof}
We prove the contraposed implication. Suppose $U$ is not a \kraj-theory. It follows that there is
an $N\ocol U\rhd {\sf S}^1_2$ such that $U^\star :=U+\verz{{\sf con}^N_n(U) \mid n \in \upomega} = U+\mho^N(U)$ is consistent.
By the Interpretation Existence Lemma, our Theorem~\ref{intexists}, we have $U^\star \rhd_{\sf res} U$.
By Theorem~\ref{sluwesmurf}, we see that $U$ is not $\imath_3$-invertible.
\end{proof}






\section{Self-Interpretations of a \kraj-Theory}\label{mosterdsmurf}
We first start with some basic facts and definitions.

\subsection{Basic Facts and Definitions}
In this section, cuts will be downwards closed under $\leq$ and closed under the arithmetical operations and the smash function.

We fix a c.e. consistent sequential theory $U$. Let $N \in \mf N_U$ iff $N \ocol U \rhd \sot$.
For $M,N\in \mf N_U$, we define:
\begin{itemize}
\item
$F\ocol M \preceq N$ expresses that $F$ is an initial homomorphic embedding of $M$ into $N$.
Here it is intended that both the domain is downward closed under $\leq^M$ and the range is downward closed
under $\leq^N$. Also it is intended that our homomorphism commutes with the smash function.
\item
$M \preceq_{\sf cut} N$ if $M$ is an $N$-cut, in other words, if $F\ocol M \preceq N$ and $F$ is the identical initial embedding on $M$.
\item
Let $M,N \in \mf N_U$,
\begin{itemize}
\item
 $F:M \preceq_U N$ iff $U \vdash F\ocol M \preceq N$. Moreover, $M \preceq_U N$ iff $F:M \preceq_U N$, for some $F$.
 \item
 $M \preceq_{{\sf cut},U} N$ iff $U \vdash  M \preceq_{\sf cut} N$.
 \item
 $F:M \prec^+_U N$ iff $F: M\preceq_U N$ and, for all $G$, $U \vdash \neg\;\; G\ocol N \preceq M$. 
 Moreover,  $M \prec^+_U N$ iff $F:M \prec^+_U N$, for some $F$.
 \end{itemize} 
 \item
 $K\ocol U\rhd U$ is \emph{initial} iff, for all $N \in \mf N_U$, there is an $M\in \mf N_U$ such that
$(K\circ M) \preceq_U N$. 
\item
$K\ocol U\rhd U$ is \emph{strongly initial} iff, for all $N \in \mf N_U$, there is an $M\in \mf N_U$ such that
$(K\circ M) \prec^+_U N$. 
\end{itemize}

We collect some easy facts.
\begin{lem}
Let $U$ be a consistent c.e. sequential theory.
\begin{enumerate}[a.]
\item
$\preceq_U$ is a partial pre-order and $\prec^+_U$ is a strict partial order.
\item
$M \preceq_U N$ and $N \prec_U^+ P$, then $M \prec_U^+ P$, and
if $M \prec^+_U N$ and $N \preceq_U P$, then $M \prec_U^+ P$. 
\item
Suppose $K\ocol U \rhd U$. Then, $ F\ocol M \preceq_U N$ implies $ F^K\ocol K \circ M\preceq_U K \circ N$.
\item
Suppose $\sigma$ is a $\Sigma^0_1$-sentence. Then $U \vdash (F\ocol M\preceq N \wedge \sigma^M) \to \sigma^N$.
\end{enumerate}
\end{lem}

A fundamental insight is as follows.

\begin{theorem}[Pudl\'ak]
Let $U$ a consistent sequential c.e. theory.
The ordering $\preceq_U$ is downward-directed, i.e., for all $M,N\in \mf N_U$, there is a $P \in \mf N_U$,
such that $P \preceq_U M$ and $P \preceq_U N$. Moreover, we can find $P$ such that 
 $P \preceq_{{\sf cut},U} M$ and $P \preceq_U N$
\end{theorem}

The proof uses the machinery of sequences to produce an isomorphic embedding of an $M$-cut in $N$.
See \cite{pudl:cuts85} or \cite{haje:meta91} or \cite{viss:what13}.

We will use the following consequence.

\begin{lem}\label{tompoessmurf}
Let $N,N',M\in \mf N_U$ and let $K:U \rhd U$. Suppose $K \circ M \preceq_U N$. Then, there is an $M'\preceq_{{\sf cut}, U}N'$, such that
$K \circ M' \preceq_U N$.
\end{lem}

\begin{proof}
We can find $M'$ with $M'\preceq_{{\sf cut},U}N'$ and $M' \preceq_U M$. Then,
$K\circ M' \preceq_U K \circ M \preceq_U N$.
\end{proof}

In Lemma~\ref{verlorenzoonsmurf},  we will need a stronger version of the Interpretation Existence Lemma.
We give the relevant definitions.
We extend the language of $U$ with a unary predicate {\sf C} and a binary predicate {\sf sat}. Let $N\in \mf N_U$. 
We will pretend that $N$ is one-dimensional leaving the modification needed for more-dimensionality and even
piece-wiseness to the reader. (Since, we are considering sequential theories, we can reduce the whole discussion to
one-dimensional anyway.)

We define:
\begin{itemize}
\item
${\sf F}^-_{U,N}$ is the theory in the language of $U$ extended with a unary predicate {\sf C} and a binary predicate {\sf sat} 
given by the following set of axioms:
\begin{itemize}
\item
$\forall x\bin {\sf C}\; {\sf form}^N(x)$, where {\sf form} expresses that $x$ is a $U$-formula.
\item
$\forall x\bin {\sf C}\, \forall y\ ({\sf sub}(y,x) \to y\in {\sf C})$, where {\sf sub} arithmetises the sub-formula relation. 
\item
${\sf C}(\gnum \phi)$, where $\phi$ is a (standard) $U$-formula and the numeral $\gnum \phi$ is derived from $N$.
\item
$\forall x\, \forall y\, ({\sf sat}(x,y) \to y\in {\sf C})$.
\item
The commutation conditions for {\sf sat} restricted to $U$-formulas $x$ in {\sf C}, where the $x$ code sequences.
\end{itemize}
\item
Suppose $M\in \mf N_V$. We have
$U \rhd_{{\sf fres},M} V$ iff $U \rhd (V+{\sf F}^-_{V,M})$.
\end{itemize} 

We have:
\begin{theorem}[Interpretation Existence Lemma 2]\label{sterkesmurf}
$\mho_U \rhd_{\sf fres} U$.
\end{theorem}

See \cite[Theorems 3.4 and 4.5]{viss:inte18}. Let $U^+:= U+{\sf F}^-_{U,N}$.

A fundamental insight is cut-intersection.
\begin{theorem}[Cut-Intersection]\label{snijsmurf}
In $U+{\sf F}^-_{U,N}$, there is an $N$-cut $\mc I$, such that, for every  $J$ defined in the $U$-language,
we have $U^+ \vdash J\preceq_{\sf cut}N \to \mc I \preceq_{\sf cut} J$.
\end{theorem}

We specify $\mc I$. Let $\varepsilon$ stand for an empty string. We define:
\[
x \in \mc I :\iff x \in \updelta_N \wedge \forall I\; \big({\sf sat}(\varepsilon,\gnum{I \preceq_{\sf cut}N}) \to
  {\sf sat}(\tupel{x},\gnum{ v_0\in I})\big)\]

\begin{rem}
We are doing things here that we did not
notationally justify. Officially one would have to do something like this. We consider application of $\lambda X. {\sf cut}_N(X)$
to $\lambda x.I(x)$. This operation can be arithmetised. Let the G\"odel number of $\lambda X. {\sf cut}_N(X) $ be $q$.
 So, we want something like ${\sf sat}(\varepsilon,\num q\cdot i)$.
 
 A second issue is that, in a sequential theory, sequences and singletons need not be unique. We can easily adapt the formula
 to take that into account.
\end{rem}

\begin{proof}
We first prove that $\mc I$ is a $U^+$-cut. We treat two cases. 
\begin{description}
\item[Downward closure]
We reason in $U^+$. Suppose $y\leq^N x$ and $x \in\mc I$. Consider any $I$ such that
${\sf sat}(\varepsilon,\gnum{I \preceq_{\sf cut}N})$. (Here we treat `$I$' as an internal variable.)
 It follows that  ${\sf sat}(\tupel{x},\gnum{v_0\in I})$. By the commutation conditions, we find
${\sf sat}(\varepsilon,\gnum{\forall u\bin I \,\forall v \bleq u\, v\in I})$.
Applying the commutation conditions again, we may conclude that ${\sf sat}(\tupel{y},\gnum{v_0\in I})$.
\item[Closure under Addition] 
We reason in $U^+$. Suppose $x$ and $y$ are in $\mc I$ and suppose
${\sf sat}(\varepsilon,\gnum{I\preceq_{\sf cut}N})$.   So,  ${\sf sat}(\tupel{x},\gnum{v_0\in I})$
and  ${\sf sat}(\tupel{y},\gnum{v_0\in I})$. We find ${\sf sat}(\varepsilon,\gnum{\forall u\bin I \,\forall v \bin I\, (u+v)\in I})$
and, thus, ${\sf sat}(\tupel{x+y},\gnum{v_0\in I})$.
\end{description}

We prove the second claim.
We reason in $U^+$. Suppose $J\preceq_{\sf cut}N$ and
 $x\in \mc I$. Then, ${\sf sat}(\tupel{x},\gnum{v_0\in J})$. Since $J$ is standard, we find $x\in J$.
\end{proof}

\subsection{Arithmetically Grounded Theories}

In this section, we briefly note some basic facts about arithmetically grounded sequential theories.
We note that these have many analogies to sequential theories with full induction. In fact, we do not know
whether they are not simply sequential theories with full induction.

\begin{ques}
Do arithmetically grounded arithmetical theories have full induction?
\end{ques}

We can derive many analogies of arithmetically grounded sequential theories with
sequential theories with full induction. Some of them are given in the next theorem.

\begin{theorem}\label{verzamelsmurf}
Suppose $W$ is a consistent, sequential, c.e., arithmetically grounded theory. Let $P$ be a witness of arithmetical groundedness. We have
\begin{enumerate}[1.]
\item
$W \vdash (\mathrm I\Delta_0 + \conOmega_1)^P$.
\item
For every $W$-sentence $\phi$ and for every $n$, we have $W \vdash {\sf prov}^P_{W, \num n}(\gnum \phi) \to \phi$.
\item
$W \vdash (\mathrm B\Sigma_1)^P$.
\item
Let $\sigma$ and $\sigma'$ be $\Sigma^0_1$-sentences and suppose ${\sf EA} \vdash \sigma\to \sigma'$. Then,
$W \vdash (\sigma \to \sigma')^P$.
\end{enumerate}
\end{theorem}

\begin{proof}
\emph{Ad \textup(1\textup):}
By a result of Wilkie, see \cite[pp 366--371]{haje:meta91}, we have $\sot \rhd \mathrm I\Delta_0$. 
So, for some $N\in \mf N_W$, we have $W \vdash  (\mathrm I\Delta_0)^N$. 
Since, $P \preceq_W N$ and  $\mathrm I\Delta_0$ can be written in $\Pi^0_1$-form, by downward persistence, we 
find that $W \vdash  (\mathrm I\Delta_0)^P$, and, thus, $W \vdash  (\mathrm I\Delta_0+\conOmega_1)^P$, since $P$ is already
$W$-provably closed under $\upomega_1$.

\medskip\noindent
\emph{Ad \textup(2\textup):}
Since $W$ is sequential, it is locally essentially reflextive. This means that, for every $W$-sentence $\phi$ and for every $n$, we can find an $N \in \mf N_W$, such that
$W \vdash {\sf prov}_{W,n}^N(\gnum \phi) \to \phi$. Since $P \preceq_W N$, we find, by the upward persistence of $\Sigma^0_1$-sentences, that
$W \vdash {\sf prov}_{W,n}^P(\gnum \phi) \to \phi$. 

\medskip\noindent
\emph{Ad \textup(3\textup):}
By (2), we have $W \vdash \mho_W^P$. By the Interpretation Existence Lemma~\ref{intexists}, we have $H:\mho_W \rhd_{\sf res} W$, where
we have a truth predicate for $H$. We have $P \preceq_W P\circ H \circ P$. We have a truth predicate $\eta$ of $H \circ P$, inside $P$,
So, we can find a liar sentence $\lambda$ such that  $U \vdash \lambda \iff \neg\, \eta(\gnum \lambda)$.
So, $U \vdash \lambda^P \iff \neg \lambda^{PHP}$. So, $P$ and $P\circ H\circ P$ are not elementarily equivalent and \emph{a fortiori} not isomorphic.
Thus, the embedding must be strictly initial. So, $P$ has a strict end-extension inside $W$.
It follows that, $W$-verifiably, we must have ${\mathrm B}\Sigma_1$ by \cite[Proposition 1]{wilk:exis89}.

\medskip\noindent
\emph{Ad \textup(4\textup):}
Suppose ${\sf EA} \vdash \sigma \to \sigma'$. By the results of \cite{wilk:sche87}, there is a $(\mathrm I\Delta_0+\conOmega_1)$-cut $J$ such that
$\mathrm I\Delta_0+\conOmega_1 \vdash \sigma^J \to \sigma'$. It follows that
$W \vdash \sigma^{JP} \to {\sigma'}^P$. So, since $P \preceq_W P \circ J$, we find $W \vdash (\sigma \to {\sigma'})^P$.  
\end{proof}

We note that Theorem~\ref{verzamelsmurf}(2) implies some {\sf PA}-lime behaviour. For example, a consistent, arithmetically grounded, sequential, c.e.
theory will not be interpretable in a finite sub-theory. Also, such a theory has no consistent finitely axiomatisable
 extension-in-the-same-language. See also \cite{viss:pean14}.

\begin{theorem}\label{hulpvaardigesmurf}
Suppose $W$ is a consistent, sequential, c.e.,  arithmetically grounded theory. Suppose $P\in \mf N_W$ is a $\preceq_W$-least element.
 Then, there is a $\Sigma^0_1$-sentence $\iota$, such that $W\nvdash \iota^P$ and $W \rhd (W+ \iota^P)$.
\end{theorem}

\begin{proof}
We find $\iota$ with $\sot \vdash \iota \iff \exists x\, {\sf prov}_{W,x}(\gnum \bot) \leq {\sf prov}_W(\gnum{\iota^P})$.
Here the witness comparison ordering applies just to the first existential quantifier of 
$\exists x\, {\sf prov}_{W,x}(\gnum \bot)$.

Our first order of business, is to show that $W \nvdash \iota^P$. Suppose 
$W \vdash \iota^P$. For some $p$, we have $W \vdash ({\sf proof}_{W}(\num p,\gnum {\iota^P}))^P$. 
It follows that $W\vdash {\sf prov}^P_{W,\num p}(\gnum \bot)$. \emph{Quod non}, by Theorem~\ref{verzamelsmurf}(2).

We define ${\sf prov}^\ast_W(y) :\iff \exists x \, {\sf prov}_{W,x}(y) \leq {\sf prov}_W(\gnum{\iota^P})$.
We note that $\iota$ is identical to ${\sf prov}_W^\ast(\gnum \bot)$ or, in other words,  ${\sf incon}^\ast(W)$.
The L\"ob conditions for ${\sf prov}^{\ast P}_W$ are $W$-verifiable. By Theorem~\ref{verzamelsmurf}(4), we
do not need to worry about the third L\"ob Condition, since we only need it to be verifiable in {\sf EA} and since ${\sf prov}^\ast_W$ is
$\Sigma^0_1$.
We have: \[W \vdash {\sf con}^{\ast P}(W) \to  {\sf con}^{\ast P}(W+{\sf incon}^{\ast P}(W)).\] 
By the Interpretation Existence Lemma, we find an interpretation $H$ with
 \[H:(W+ {\sf con}^{\ast P}(W))\rhd (W+ {\sf incon}^{\ast P}(W)).\]
Let $K:={\sf ID} \tupel{{\sf incon}^\ast(W)}H$, where {\sf ID} is the identity interpretation on $W+{\sf con}^{\ast P}(W)$.
So, $K:W \rhd (W+{\sf incon}^{\ast P}(W))$, i.e., $K:W \rhd (W+\iota^P)$.
\end{proof}

\begin{cor}\label{smidsmurf}
Suppose $W$ is a consistent, sequential, arithmetically grounded theory. Then, there is a non-initial self-interpretation $K^\ast$ of $W$.
\end{cor}

\begin{proof}
Let $P$ witness arithmetical groundedness.
Let $K^\ast:W \to W$ be based on the translation of a witnessing interpretation of  $W \rhd (W+ \iota^P)$.
Suppose we had $M \in \mf N_W$ with $K^\ast \circ M \preceq P$. It follows that $K^\ast \circ P \preceq P$.
Since, according to $W$, we have $\iota$ inside $K^\ast \circ P$, by upwards persistence, it follows that
$W \vdash \iota^P$. \emph{Quod non.} 
\end{proof}

\begin{rem}
Theorem~\ref{hulpvaardigesmurf} is based on a combination of ideas. The first is Feferman's Theorem of the interpretability of inconsistency.
This is combined with a variation of a familiar trick to produce a provability predicate for extensions of {\sf PA} where the unprovability
of inconsistency is guaranteed. This trick in its turn uses a variant of the Friedman-Goldfarb-Harrington
argument. The true priority for this argument goes to John Shepherdson in \cite{shep:rep61}. 

The alternative provability predicate produced
in the case of extensions of {\sf PA} is Fefermanian. This means that it can be viewed as ordinary provability with an alternative representation of
the axiom set. The predicate we use here is not Fefermanian because we also employ restriction of the complexity of the formulas used in the proofs.

Arithmetical groundedness is a heavy demand. The reader may well wonder whether we cannot do better for a Feferman-FGH-style result. Indeed, with
a bit more delicate work, we can prove a similar result for reflexive, c.e., sequential theories. See Appendix~\ref{nademaaltijd}. 
\end{rem}

\begin{rem}
Since we have Theorem~\ref{verzamelsmurf}\textup(4\textup), we can verify Solovay's proof for the completeness of provability logic for
${\sf prov}_U^{\ast N}$ for interpretations in $U$ using \cite{dejo:proo91}.
\end{rem}

\subsection{The Characterisation}
We will  prove the following characterisation of consistent c.e. \kraj-theories.

\begin{theorem}\label{grotesmurf}
Let $U$ be a consistent c.e. sequential theory. The following are equivalent:
\begin{enumerate}[a.]
\item
 $U$ is a \kraj-theory.
 \item
All self-interpretations of $U$ are strongly initial.
\item
All self-interpretations of $U$ are initial.
\end{enumerate}
\end{theorem}

We note that the definition of a \kraj-theory is based on a great number of very specific design choices.
We have to choose, at the object-level, an appropriate syntax, a proof system, the desired notion of formula complexity,
the arithmetisation. In contrast, the only conventions required for the characterisation are at the meta-level. These are the kind of choices shared
by all logic. So, the characterisation goes a long way towards a `coordinate-free' way of understanding being a \kraj-theory.

We will prove a number of theorems that can be combined to a proof of Theorem~\ref{grotesmurf}.
The proof below shows how Theorem~\ref{grotesmurf} follows from these later theorems.

\begin{proof}[Proof of Theorem~\ref{grotesmurf}]
(a) $\To$ (b): This is Lemma~\ref{neergaandesmurf}.

\medskip\noindent
(b) $\To$ (c):  This is trivial.

\medskip\noindent
(c) $\To$ (a): Suppose $U$ is a consistent sequential c.e. theory and suppose all self-interpretations of $U$ are initial.
Suppose, in order to obtain a contradiction, that
$U$ is not a \kraj-theory. This tells us that, for some $N \in \mf N_U$, the theory $U+\verz{{\sf con}^N_n(U) \mid n\in \upomega}$ is consistent.
By Theorem~\ref{slimmesmurf}, there is a $\zeta$ such that $U+\zeta$ is consistent and $U + \zeta \vdash \verz{{\sf con}^N_n(U) \mid n\in \upomega}$.
Theorem~\ref{oneupsmurf} tells us that all self-interpretations of $U+\zeta$ are initial. Moreover, by Lemma~\ref{verlorenzoonsmurf},
$U+\zeta$ is arithmetically grounded. But his contradicts the consistency of $U+\zeta$ in the light of Corollary~\ref{smidsmurf}.
\end{proof}

We prove Lemma~\ref{neergaandesmurf}.

\begin{lem}\label{neergaandesmurf}
Let $U$ be a consistent c.e. sequential \kraj-theory. Then any self-interpretation $K$ of $U$ is
strongly initial.
\end{lem}

\begin{proof}
We work in $U$. Consider any $N\in \mf N_U$. Let $P \preceq_U N$ and $P\preceq_U K\circ N$.
Let $G$ witness  $P\preceq_U K\circ N$.
We know that $U \vdash {\sf incon}^P_n(U)$, for some $n$.
Let $M$ in $\mf N_U$ be an $N$-cut such that $U \vdash {\sf con}^M_n(U)$. It follows that $U \vdash  {\sf con}^{MK}_n(U)$, since we have
$U$ inside $K$.

We want to show that $K\circ M \preceq_U P \preceq_U N$. We reason in $U$. It suffices that every element of $K\circ M$ is in the range of $G$.
Let $a$ be an arbitrary element of $K\circ M$ and let  $b$ be a witness of ${\sf incon}^P_n(U)$. It follows that $G(b)$ is a witness of ${\sf incon}^{NK}_n(U)$.
Since $ {\sf con}^{MK}_n(U)$, we find that $a<^{NK}G(b)$. Ergo, $a$ is in the range of $G$.

Let $ F: P \preceq_U N$.
We reason in $U$.
Suppose $G: N\preceq_U K\circ M$. If follows that $G\circ F: P \preceq K \circ M$.   We have ${\sf incon}^P_n(U)$, so ${\sf incon}^{MK}_n(U)$.
\emph{Quod non.} 
\end{proof}

We prove that finite extensions preserve the property that all self-interpretations are initial.

\begin{theorem}\label{oneupsmurf}
Let $U$ be a consistent c.e. sequential theory and let $\phi$ be any $U$-sentence such that $U+\phi$ is consistent.
Suppose all self-interpretations of $U$ are initial.
Then, all self-interpretations of $U+\phi$ are initial.
\end{theorem}

\begin{proof}
Let $N_0\in \mf N_U$ with indwelling translation $\nu_0$. 
Let $K$ be a self-interpretation of $U+\phi$ and let $N \in \mf N_{U+\phi}$. 
Let $\kappa$ be the indwelling translation of $K$ and let $\nu$ be the indwelling translation of $N$.
Let $K'$ be the self-interpretation of $U$ based on $\kappa \tupel{\phi} {\sf id}$ and let $N'$ be the element of
$\mf N_U$ based on $\nu \tupel{\phi}\nu_0$. We can find an $N^\circ\in \mf N_U$ such that  $ F: K'\circ N^\circ \preceq_U N'$. 
 Let $N^\ast := {\sf E}_{U(U+\phi)} \circ N^\circ  \in \mf N_{U+\phi}$.
Then, $U+\phi \vdash F: K\circ N^\ast \preceq N$.
\end{proof}

We prove the promised Lemma~\ref{verlorenzoonsmurf}.

\begin{lem}\label{verlorenzoonsmurf}
Let $U$ be a consistent c.e. sequential theory. Suppose all self-interpretations of $U$ are initial. Suppose further that
$U+\zeta$ is consistent and $U+\zeta \vdash
\verz{{\sf con}^N_n(U) \mid n \in \upomega}$, for some $N\in \mf N_U$. Then, $U+\zeta$ is arithmetically grounded.
\end{lem}

\begin{proof}
We assume the conditions of the theorem.
Let $K^+:(U+\zeta) \rhd (U+{\sf F}^-_{U,N})$ and let $K$ be the reduct of $K^+$ to the $U$-language.
We take $K^\ast := K\tupel{\zeta} {\sf E}_{U(U+\neg\, \zeta)}:U \rhd U$.

We find that $K^\ast$ is initial. 
Consider any $N'\in \mf N_{U+\zeta}$. Let $N^\ast := N'\tupel{\zeta} N'': U \rhd \sot$, where
 $N'':= {\sf E}_{U(U+\neg\,\zeta)}\circ N \in \mf N_{U+\neg\, \zeta}$.
By Lemma~\ref{tompoessmurf}, it follows that there is an $M\preceq_{{\sf cut},U} N$, 
such that $ F:K^\ast \circ M \preceq_U N^\ast$, for some $F$. 
It follows that $ F: K^+\circ M' \preceq_{U+\zeta} N'$, where $M'= {\sf E}_{U(U+\zeta)}\circ M\in \mf N_{U+\zeta}$.
Since $M'$ is in the $U$-language, it follows that (\dag) $  F: K\circ M' \preceq_{U+\zeta} N'$,

Theorem~\ref{snijsmurf}, gives us $U+{\sf F}^-_{U,N} \vdash M \preceq_{\sf cut} N \to \mc I \preceq_{\sf cut} M'$,
so \[(\ddag) \;\;\;\;   K^+\circ \mc I \preceq_{{\sf cut},(U+\zeta)} K \circ M'.\]
Combining, (\dag) with (\ddag), we find that $K^+\circ \mc I \preceq_{U+\zeta} N'$. So, $K^+ \circ \mc I$ is a desired least element of $\preceq_{U+\zeta}$.
\end{proof}


\appendix
\section{A Feferman-FGH-style Result}\label{nademaaltijd}

Axiom sets of theories are supposed to be represented by a $\Delta_1^{\sf b}$ formula over $\sot$.
We prove the following.

\begin{theorem}
Let $W$ be c.e. and sequential. Let $N \in \mf N_W$. Suppose $W$ is $N$-reflexive, i.e.,
$W\vdash {\sf con}^N_n(W)$, for each $n$. Then, there is a $\exists\Sigma^{\sf b}_1$-sentence
$\iota$, such that $W \nvdash \iota^N$ and $W \rhd (W+ \verz{\iota^M\mid M\in \mf N_W})$.
\end{theorem}

\begin{proof}
We find $\iota$ with \[\sot \vdash \iota \iff \exists x\, ({\sf prov}_{W,x}(\gnum \bot) \wedge \forall y\bles |x|\, \neg \, {\sf proof}_W(y,\gnum {\iota^N})).\]
Here $\forall y\bles |x|$ is sharply bounded. 
We note that we can take $\iota$ to be $\exists\Sigma_1^{\sf b}$.

We show that $W \nvdash \iota^N$. Suppose 
$W \vdash \iota^N$. For some $p$, we have \[W \vdash {\sf proof}^N_{W}(\num p,\gnum {\iota^N}).\]
We note that, in $W$, it follows from $\iota^N$ and ${\sf proof}^N_{W}(\num p,\gnum {\iota^N})$, that the promised witness $x$ of $\iota^N$
must satisfy $|x| \leq \num p$ and, thus, $x \leq 2^{\num p}$.
Let, in the real world, $q := 2^{p}$.
So, we will certainly have  $W\vdash {\sf prov}^N_{W,\num q}(\gnum \bot)$. \emph{Quod non}.

Let $M \in \mf N_W$. We define 
\begin{itemize}
\item
${\sf prov}^\ast_W(z) :\iff \exists x \, ({\sf prov}_{W,x}(z) \wedge \forall y\bles |x| \neg\, {\sf proof}_W(y,\gnum{\iota^N}))$.
\item
${\sf prov}_W^\circ(\gnum\phi) :\iff {\sf prov}_W^\ast (\gnum{\phi^M})$, where $\phi$ is arithmetical. 
\end{itemize}
We note that $\iota$ is ${\sf incon}^\ast(W)$. To increase readability, we have hidden the dependence of ${\sf prov}_W^\circ$ on $M$. 

We can easily verify the L\"ob Conditions for ${\sf prov}_W^\circ$.
We find, for sufficiently large $n$, that:{\small
\begin{eqnarray*}
\sot \vdash {\sf prov}_{W,n}(\gnum{ {\sf con}^{\ast M}(W)}) & \to & 
\big ({\sf prov}_{W,n}(\gnum {{\sf con}^{\ast M}(W)}) \wedge {\sf prov}^\circ_{W} (\gnum{{\sf con}^\circ(W)}) \big)\\
&\to & \big({\sf prov}_{W,n} (\gnum{{\sf con}^{\ast M}(W)}) \wedge{\sf incon}^\circ (W) \big)\\
&\to & \big({\sf prov}_{W,n} (\gnum{{\sf con}^{\ast M}(W)}) \wedge{\sf incon}^\ast (W) \big)\\
&\to & \big({\sf prov}_{W,n} (\gnum{{\sf con}^{\ast M}(W)}) \wedge {\sf prov}_{W,n}(\gnum{{\sf incon}^\ast (W)}) \big)\\
&\to & {\sf prov}_{W,n} (\gnum \bot)
\end{eqnarray*}
}
We may conclude that $\sot \vdash {\sf con}_n(W) \to {\sf con}_n(W+\iota^M)$. 
So,
$W \vdash {\sf con}^N_n(W+\iota^M)$. 
We note that at this point, we may drop the assumption that $n$ be sufficiently large.

Consider an effective enumeration $M_0,M_1,\dots$ of $\mf N_W$.
Since, for each $k$, we can find an $M$ such that $M \preceq_W M_i$, for all $i<k$, we find that
\[W \vdash {\sf con}^N_n(W+\iota^{M_0} + \cdots + \iota^{M_{k-1}}).\]

By the Interpretation Existence Lemma, we find 
 $W\rhd (W+ \verz{\iota^{M_i} \mid i \in \upomega})$.
\end{proof}

\end{document}